\documentclass[review]{elsarticle}

\usepackage{lineno,hyperref}
\usepackage{amsmath}
\usepackage{amsfonts}
\usepackage{threeparttable}
\usepackage{multirow}
\usepackage{booktabs,color}
\usepackage{algorithm,algpseudocode} 
\modulolinenumbers[5]
\usepackage{caption}
\usepackage{subfigure}
\newtheorem{definition}{Definition}

\journal{Journal of \LaTeX\ Templates}

\begin{document}

\begin{frontmatter}

\title{Learning to Select Cuts for Efficient Mixed-Integer Programming}

\author[mymainaddress]{Zeren Huang}
\author[mymainaddress]{Kerong Wang}

\author[mythirdaddress]{Furui Liu\corref{cor1}}
\ead{liufurui2@huawei.com}
\author[mythirdaddress]{Hui-Ling Zhen\corref{cor1}}
\ead{zhenhuiling2@huawei.com}
\author[mymainaddress]{Weinan Zhang\corref{cor1}}
\ead{wnzhang@sjtu.edu.cn}
\author[mythirdaddress]{Mingxuan Yuan}
\author[mythirdaddress]{Jianye Hao}
\author[mymainaddress]{Yong Yu}
\author[mysecondaryaddress]{Jun Wang}
\cortext[cor1]{Corresponding author}
\address[mymainaddress]{Shanghai Jiao Tong University}
\address[mysecondaryaddress]{University College London}
\address[mythirdaddress]{Noah's Ark Lab, Huawei Technologies}

\begin{abstract}
Cutting plane methods play a significant role in modern solvers for tackling mixed-integer programming (MIP) problems. Proper selection of cuts would remove infeasible solutions in the early stage, thus largely reducing the computational burden without hurting the solution accuracy. However, the major cut selection approaches heavily rely on heuristics, which strongly depend on the specific problem at hand and thus limit their generalization capability. In this paper, we propose a data-driven and generalizable cut selection approach, named \textsc{Cut Ranking}, in the settings of multiple instance learning. To measure the quality of the candidate cuts, a scoring function, which takes the instance-specific cut features as inputs, is trained and applied in cut 
ranking and selection. In order to evaluate our method, we conduct extensive experiments on both synthetic datasets and real-world datasets. Compared with commonly used heuristics for cut selection, the learning-based policy has shown to be more effective, and is capable of generalizing over multiple problems with different properties. \textsc{Cut Ranking} has been deployed in an industrial solver for large-scale MIPs. In the online A/B testing of the product planning problems with more than $10^7$ variables and constraints daily, \textsc{Cut Ranking} has achieved the average speedup ratio of 12.42\% over the production solver without any accuracy loss of solution. 
\end{abstract}



\end{frontmatter}


\section{Introduction}
Combinatorial optimization (CO) is a subclass of optimization problems, where the goal is to find the optimal solution with respect to a given objective function from a finite candidate solution set. Due to its combinatorial nature (for example, integer constraints), it is usually NP hard and can mostly be formulated as mixed-integer programming (MIP) problems \cite{bixby2012brief,richards2005mixed,achterberg2013mixed,bixby2004mixed}. It covers a wide range of industry applications such as production planning, scheduling and manufacturing \cite{wu2013mixed, schouwenaars2001mixed, keha2009mixed, amaral2019mixed}.


The difficulty of solving MIP problems lies on the non-convexity of its feasible region, which makes the general MIPs unsolvable in polynomial time. Instead of solving the MIP directly, one usually solves the corresponding LP relaxations first, and then performs rounding to generate the approximately optimal solution \cite{wolsey2007mixed}. To facilitate such a process, a classic approach is the cutting plane method, which generates valid inequalities in the LP iterations to cut off the fractional solutions or the infeasible integer solutions, so that the convergence to optimum is accelerated \cite{marchand2002cutting}. Another approach to solve MIP is the branch-and-bound algorithm \cite{lawler1966branch, ris2010u}, which creates branches by selecting variables to add rounding bounds to form two sub-LP problems (for a integer variable $x_i$ with fractional value $v$ as example, add two rounding bounds $x_i \geq \lceil v \rceil$ and $x_i \leq \lfloor v \rfloor$), and then solve these sub-LP problems. In modern MIP solvers, the cutting plane technique is often combined with the branch-and-bound method to constitute the branch-and-cut framework \cite{mitchell2002branch}, where each branch contains cutting planes as additional constraints. 
Heuristics are also employed in this process to cope with problems such as branching variable selection and cutting plane (cut) selection, i.e., selecting the most promising cuts to add. Such heuristics in MIP solvers are usually manually designed and heavily dependent on the problem. As a result, it shows high vulnerability with respect to the structure or the size of the MIP \cite{achterberg2012rounding}.

As a promising methodology for address above issues, many recent works \cite{balcan2018learning,khalil2016learning,gasse2019exact, tang2020reinforcement} have leveraged machine learning (ML) techniques to construct efficient heuristics which is problem independent, and the majority are focusing on decision problems in the branch-and-bound algorithm, thus leaving room for machine learning methods on the direction of cut selection.

A good set of cuts is essential for the efficiency of the CO algorithms.  Cuts serve the purpose to reduce the LP solution space, which results in a smaller tree in branch-and-cut algorithm so that the number of nodes to be searched are significantly reduced. However, excessive quantity of cuts causes heavy computational workload on solving corresponding LP problems. As a consequence, deriving a good cut selection policy is of high value to the community, which has unfortunately received few attention so far. This motivates us to develop a general data-driven, machine learning based cut selection algorithm.

The purpose of this work is to construct an efficient and generalizable cut selection policy based on machine learning. The basic idea is to learn a scoring function that can measure the quality of cuts, and we formulate it as a \emph{cut ranking} problem, in which we score each generated cut by a learned scoring function, and select a subset of cuts with the highest scores. However, such a task is non-trivial, with several remarkable technical challenges. 
First, in many cases, labels for individual cuts are not easy to obtain, since the impact of a single cut on MIP is relatively weak and imperceptible. Labeling good cuts individually may be infeasible. Thus, the cut selection problem naturally fits the scenarios of \emph{multiple instance learning} (MIL), in which the collection of labels is at the bag level \cite{carbonneau2018multiple,foulds2010review,cheplygina2015multiple}.
The training instances are organized into sets (also called bags), and the label is not assigned to any individual instance, but to the bag of cuts to measure the overall quality.
Another important problem is the generalization ability of the learned scoring function. To enable the scoring function module to generalize to new problems, we design both static and dynamic problem-specific cut features as the inputs. Inspired by the training process of Reinforcement Learning (RL), we collect the supervised labels in an exploratory way \cite{arulkumaran2017brief}.
Our proposed method is named as \textsc{Cut Ranking}.

To summarize, the technical contributions of our work are threefold.
\begin{enumerate}
	\item We propose a novel \textsc{Cut Ranking} method for cut selection in the settings of multiple instance learning, which is suitable for the nature of cut selection tasks.
	\item We study generalization ability to the cut selection policy since the designed cut features are determined by the characteristics of MIPs.
	\item The \textsc{Cut Ranking} module can be applied as a subroutine in the branch-and-cut algorithm, which is generally adopted by modern MIP solvers.
\end{enumerate}

The extensive experiments on various MIP problems demonstrate the superiority of \textsc{Cut Ranking} over previous solutions in terms of solving time and the node size of the branch-and-bound tree.
Furthermore, we deploy \textsc{Cut Ranking} on Huawei's proprietary industrial large-scale MIP solver for production planning problems with more than $10^7$ variables and constraints daily. The A/B testings show that our solution can reduce the overall solving time by 14.98\% and 12.42\% on average for offline and online phases respectively without the accuracy loss of solution, which is a significant acceleration for the industrial solver.

\section{Related works}
The traditional approaches to tackle the MIPs mainly include: branch-and-bound \cite{lawler1966branch,clausen1999branch,morrison2016branch} and cutting-plane methods \cite{marchand2002cutting, cornuejols2008valid}. They are widely deployed in modern solvers as the core algorithm for solving problems. However, in the era of big data, large scale problems with a lot of variables   are often encountered, and those approaches suffer from very low efficiency.   For scalability and speeding up the solvers, they are usually enhanced with  heuristics, which are often designed by experts, based on the unique property of the problems at hand, and are not transferable or reusable when one switches to a new situation.

Therefore, there emerges a need for generalizable methods that are ubiquitous applicable to MIPs. Attentions are thus paid to machine learning and other data driven science due to their generalizability.
Given  training data, intelligent models are able to learn to solve the problem, with good performance on unseen data.  Based on the role that the machine learning model takes, related literatures can be categorized into two clusters \cite{bengio2021machine}. The first cluster contains methods that use ML models to replace traditional solving techniques. They can directly solve the MIP problems (such as TSPs), in which an end-to-end learning model is often used to predict the solution given the problem instance. Vinyals et al.~\cite{vinyals2015pointer} proposed the pointer network with a sequence-to-sequence architecture, and train the model through supervised learning. Bello et al.~\cite{bello2016neural} introduced a reinforcement learning method to train the pointer network. Other literatures including \cite{kool2018attention,nazari2018reinforcement} also use the sequence-to-sequence architecture to tackle the vehicle routing problems. More recently, Nair et al.~\cite{nair2020solving} showed that the deep learning model is able to predict a good partial solution for MIP problems. 
The aforementioned works often use blackbox and unexplainable AI models, supported by empirical evidences but no theoretical guarantees. As a result, they sometimes show performance vulnerability in solving MIPs.

For the second cluster, in which our work can be placed, the main algorithmic framework is based on the traditional optimization algorithm, and machine learning is used to improve the  heuristics. During the solution process, the ML model is repeatedly called to assist in making decisions. There have been multiple studies \cite{balcan2018learning,khalil2016learning,gasse2019exact} about learning a branching policy in the context of branch-and-bound, in which the policy is usually trained through imitation learning to approximate a powerful heuristics named strong branching, which is effective but too slow. Khalil et al.~\cite{khalil2016learning} addressed the branching variable selection problem as a ranking problem. 

In addition to learning a branching strategy, there are also studies on learning other core elements in the branch-and-bound framework. He et al.~\cite{he2014learning} proposed to learn a node selection policy to improve the heuristics. Khalil et al.~\cite{khalil2017learning} adopted a ML-based model to decide running a given heuristic or not at a branching node.  

To our knowledge, the direction of deriving a ML-based cut selection policy for MIP problems have been rarely explored except Tang et al.~\cite{tang2020reinforcement}, in which the authors introduced a MDP formulation for the problem of iteratively selecting cuts for a MIP, and train a reinforcement learning (RL) agent using evolutionary strategies. Our work differs from it in two aspects: first, we propose a ranking formulation, and model the learning problem in multiple instance learning settings; second, the main goal of Tang et al.~\cite{tang2020reinforcement} is to improve the efficiency of the cut selection module, i.e. to reduce the total number of cuts added, while the objective of our work is to improve the performance of the optimization algorithm, i.e. to reduce the total running time.  

\section{Background}

Our algorithm combines MIL technique with the branch-and-cut framework. In this section, we  firstly introduce the background of MIP and branch-and-cut framework. Then, MIL related techniques are presented.
\subsection{Mixed-Integer Programming Background}
\paragraph{\textbf{Mixed-Integer Programming Problem}} 

The general formulation of a MIP is as follows:
\begin{equation}
\label{mip_formulation}
	\underset{\mathbf{x}}{\arg \min}\left\{\mathbf{z}^{\top} \mathbf{x} \mid \mathbf{A} \mathbf{x} \leq \mathbf{b},  \mathbf{x} \in \mathbb{Z}^{p} \times \mathbb{R}^{n-p}\right\},
\end{equation} 
where $\mathbf{x}$ is the vector of decision variables, $\mathbf{z} \in \mathbb{R}^n$ is the objective coefficient vector, $\mathbf{b} \in \mathbb{R}^{m}$ is the right-hand side vector, $\mathbf{A} \in \mathbb{R}^{m \times n}$ is the constraint matrix.

\paragraph{\textbf{Branch-and-Cut}}

The branch-and-cut (BC), a combination of two classical algorithms, is widely adopted in modern MIP solvers. To address the difficulty brought by the nonconvexity  on searching the optimal solution,    it builds a  search tree with each node corresponding to a linear programming problem. At the beginning of the algorithm, the root node corresponds to the problem with all integer constraints dropped. Then, it iteratively generates child nodes by selecting variables to branch on, that is, adding new constraints (bounds) on it. Along the paths of the tree, the space of solution is regularized by more constraints. For the problems corresponding to each node, valid cuts are added to assist searching. 
Taking the root node for example, assume that the added cut set $C'=\{\alpha_i^\top x \leq \beta_i\}_{i=1}^{|C'|}$ is a subset of generated cut set $C$, the optimization problem becomes
\begin{equation}
	\label{mip_formulation_2}
	\underset{\mathbf{x}}{\arg \min}\left\{\mathbf{z}^{\top} \mathbf{x} \mid \mathbf{A} \mathbf{x} \leq \mathbf{b}, 
	\mathbf{\alpha}^\top \mathbf{x} \leq \mathbf{\beta},
	 \mathbf{x} \in \mathbb{Z}^{p} \times \mathbb{R}^{n-p}\right\}.
\end{equation} 

The algorithm terminates when there exists no feasible solutions for one node, or we cannot obtain a better solution than the optimal one found so far. It is worth noting that the cuts are primarily added at the root node, which will often bring significant improvements.

\paragraph{\textbf{Metric for Cut Quality}}
For a given MIP, its solvability is defined as the capability of being solved, which relates to the size, the structure and other problem properties. 
A formal quantitative description of problem solvability is provided in Definition $\ref{def:solvability}$.

\begin{definition}[Problem Solvability]
	\label{def:solvability}
	Let $O$ be the set of all the feasible optimization algorithms for solving MIPs, and $S$ be the set of all the feasible MIP solvers. Assume that the computing environment is kept the same, we define the problem solvability of a MIP problem $\chi$ with respect to parameters $A,b,z$ and integrality constraints as
	\begin{equation}
		PS(\chi) = \mathbb{E}_{o\sim O,s \sim S} \Big[\frac{1}{T_{\chi}}\Big|o,s\Big],
	\end{equation}
	where $T_{\chi}$ is the solving time of the MIP problem $\chi$.
\end{definition}

To measure the quality of selected cut subset $C'$, we propose a metric named  \emph{problem solvability improvement} (PSI), which is calculated after adding cuts via
\begin{equation}
\label{equ:PSI}
	PSI  = PS(\chi') - PS(\chi),
\end{equation} 
where $\chi'$ and $\chi$ represent the MIP with and without cuts, respectively.

In real practice, it is impractical to calculate $PSI$ since obtaining the problem solvability is infeasible. However, when the optimization algorithm, the MIP solver and the computing environment are fixed, we can substitute $PSI$ with the reduction ratio of solution time as the feedback $r$ of selecting $C'$:
\begin{equation}
\label{equ:RRST}
	r = \mathbb{E}\Big[\frac{T_{\chi} - T_{\chi'}}{T_{\chi}}\Big|o, s\Big] ,
\end{equation} 
where $o$ is the optimization algorithm, $s$ is the solver. A higher value of $r$ implies a higher-quality cut subset for the MIP.

\paragraph{\textbf{Typical Cut Types}}
In cutting plane tasks, there exist various types of cuts which can be generated. Here, we list several typical types of cuts:
\begin{itemize}
	\item \textbf{Cover Cut.} For a set of binary variables $X=\{x_1, x_2, \ldots, x_k\}$, a so-called knapsack constraint takes the form as 
	\begin{equation}
		a_1 x_1 + a_2 x_2 + \cdots + a_k x_k \leq b,
	\end{equation}
	where $a_1, a_2, \ldots, a_k, b$ are all non-negative. Let $X' = \{x_1',x_2',\ldots,x_l'\} \subset X$. A minimal cover cut related with the above knapsack constraint is of the form as 
	    \begin{equation}
		x'_1 + x'_2 +\cdots + x'_l \leq l-1.
	\end{equation}
	\item \textbf{Gomory Cut.} The gomory cuts are generated from the rows of the simplex tableau, returned by the simplex algorithm for solving LPs. Here we use the similar notations as in Tang et al. \cite{tang2020reinforcement}. Denote the constraint matrix and the constraint vector of the tableau as $A^{'}$ and $b^{'}$, respectively. For the $i$th row, the corresponding gomory cut can be generated by applying integer rounding as 
	\begin{equation}
		(-A'_i + \lfloor A'_i \rfloor) x \leq -b'_i + \lfloor b'_i \rfloor .
	\end{equation}
	
	\item \textbf{Clique Cut.} For a set of binary variables $X=\{x_1, x_2, \ldots, x_k\}$, a clique cut is of the form as
	\begin{equation}
		x_1 +x_2 + \cdots + x_k \leq 1,
	\end{equation}
	where at most one variable can be positive.
	
\end{itemize}

For a more detailed introduction of other cut types, one can refer to the surveys for cutting planes \cite{marchand2002cutting, cornuejols2008valid}. Note that these different types of cuts are enabled in general MIP solvers. For a given MIP instance, we can use the cut generators incorporated in the MIP solver to generate the candidate cuts.

\subsection{Multiple Instance Learning}

Multiple instance learning (MIL) concerns the problem of supervised learning where the model prediction and training are put at the level of bag of instances \cite{carbonneau2018multiple}. Each bag is composed of multiple unlabeled training instances. The goal is to predict the labels of unseen data at the bag level or at the instance level.

For the binary classification problems, where the label is positive or negative, the standard setting of MIL is that, bags containing at least one positive instances are assigned positive labels while bags containing only negative instances are assigned with negative labels \cite{carbonneau2018multiple}. This can be relaxed to the collective assumption, which is related to problems where the label assignment is determined by more than one instances.
The MIL measures the effect of a set of instances by interpreting labels of bags,  and this naturally fits the scenario of our cut selection. For large-scale problems, the effect of a single cut to the solution is rather minor and imperceptible. Therefore, we use MIL related techniques as our ML model for cut selection.

\section{Methodology: Cut Ranking}

\subsection{A Cut Ranking Formulation in MIL Settings}
In the branch-and-cut framework, we introduce cut selection at the root node of the branching tree. Since the majority of cuts are added to the root LP in general cases, the cuts are disabled at other sub-nodes for a better evaluation of effects on the algorithm. However, we argue that the learned cut selection policy can generalize to other sub-nodes (sub-MIPs) due to that the designed cut features are of the same properties with a MIP instance. 

\begin{definition}[MIP Problem Property] 
	\label{def:mip_property}
	For a MIP with parameters $M=\{A, b, z\}$, after its root LP relaxation being solved, the LP solution $x_{LP}^{*}$ is accessible. We define the problem property of the MIP before cut selection as $P = \{ x_{LP}^{*}, M \}$.
\end{definition}

For a given MIP, its problem property $P$ is defined in Definition \ref{def:mip_property}. Let $C=\{c_1,c_2,\ldots, c_l\}$ be the candidate cut set generated at the root node, the cut selection problem is equivalent to select an optimal cut subset $C^{*}$ with respect to the problem solvability improvement (PSI) mentioned in Equation \ref{equ:PSI}:
\begin{equation}
	C^* = \underset{C'}{\arg \max }\left\{PSI \mid C' \subset C, P\right\}.
\end{equation}
Due to its combinatorial structure, finding an exact solution is intractable, especially when the size of $C$ becomes larger. 

To tackle such a problem, we present a \textsc{Cut Ranking} formulation in the branch-and-cut framework. The process of cut ranking involves the training phase and the test phase. In the training phase, the learning process is modeled in MIL settings, that is, the training data are grouped into bags, and the label assignment is at the bag level. Specifically, each bag consists of several cuts sampled from the candidate cut set, and the bags are not disjoint. Denote $u_i \in \mathbb{R}^{H}$ as the feature vector of cut $c_i$, which can be derived given the problem property $P$ and the cut parameters. Let $X=\{x_1, x_2, \ldots, x_h\}$ be the set of feature vectors for the collected bags, where each bag feature vector $x_i \in \mathbb{R}^{H}$ can be constructed through a feature mapping function $\phi(\cdot)$, taking the aggregated features of cuts within the bag as the input. 
Let $Y=\{y_1, y_2, \ldots, y_h\}$ be the labels of $X$. Given the training set, our goal is to train a scoring function $f_{\theta}(u)$ which can predict the score for each candidate cut $u$, with the cut feature vector as the input. 

In the test phase, for the given MIP, we use the trained scoring function $f_{\theta}$ to assign scores to all the generated candidate cuts, and select the top $K\%$ cuts with the highest scores (\emph{K} is a hyper-parameter).

\subsection{Constructing Training Data}
The training data is composed of features $X$ and labels $Y$ at the bag level. To construct it, we first collect the training samples (each sample corresponds to a bag of cuts) using a certain searching strategy on multiple randomly generated instances; next we extract the designed instance-specific cut features for the cuts within each training sample; after that, we construct the bag features using the aggregated features of cuts in each bag; finally, we assign binary label to each sample through a designed labeling scheme.

\subsubsection{Strategies of Collecting Training Samples}
\label{sec:strategy}
 For a given MIP, after its root LP being solved, a set of candidate cuts $C$ are generated by the cut generators incorporated in the MIP solver first. Note that for a MIP instance, when the MIP solver is fixed, the generated candidate cuts are also fixed. Denote hyper-parameter $K\%$ as the cut selection ratio, a subset $C'$ is selected using a stochastic cut selection policy, and under the selection ratio. After adding the selected cuts, the algorithm continues until terminating and returns the solution time. To measure the quality of the selected subset $C'$, since the pre-defined problem solvability is unavailable in practice, we use the reduction ratio of solution time mentioned in Equation \ref{equ:RRST}, as an alternative to $PSI$. The reduction ratio of solution time can also be regarded as the feedback $r$ of selecting a cut subset.

 For each MIP instance, we repeat running the solver for multiple times, and collect a number of training samples. Since our cut selection policy is stochastic, we are able to explore different cutting results and obtain training samples with much diversity. Note that though collecting training samples leads to multiple rounds of execution of the MIP solver, the whole process is conducted in an offline way, and thus the incurred training cost is acceptable. The collected training sample can be seen as a tuple $(P, C', r)$ consisting of the MIP property $P$, the selected cut subset $C'$ (bag of cuts), and the feedback $r$.  To improve the exploration and also the data efficiency, we collect the training samples based on two strategies, random sampling and active sampling, which are similar to the searching strategies adopted by Bello et al. \cite{bello2016neural}.

 \paragraph{\textbf{Random Sampling}} 
 The cut selection is based on a fixed stochastic policy, which randomly selects a subset of cuts to add to the MIP. For each generated MIP training instance, the algorithm is repeatedly called for a certain number of times, in which we apply random sampling to collect the initial training samples.
 
\begin{figure}[htb]\begin{center}
\includegraphics[width=1.0\textwidth]{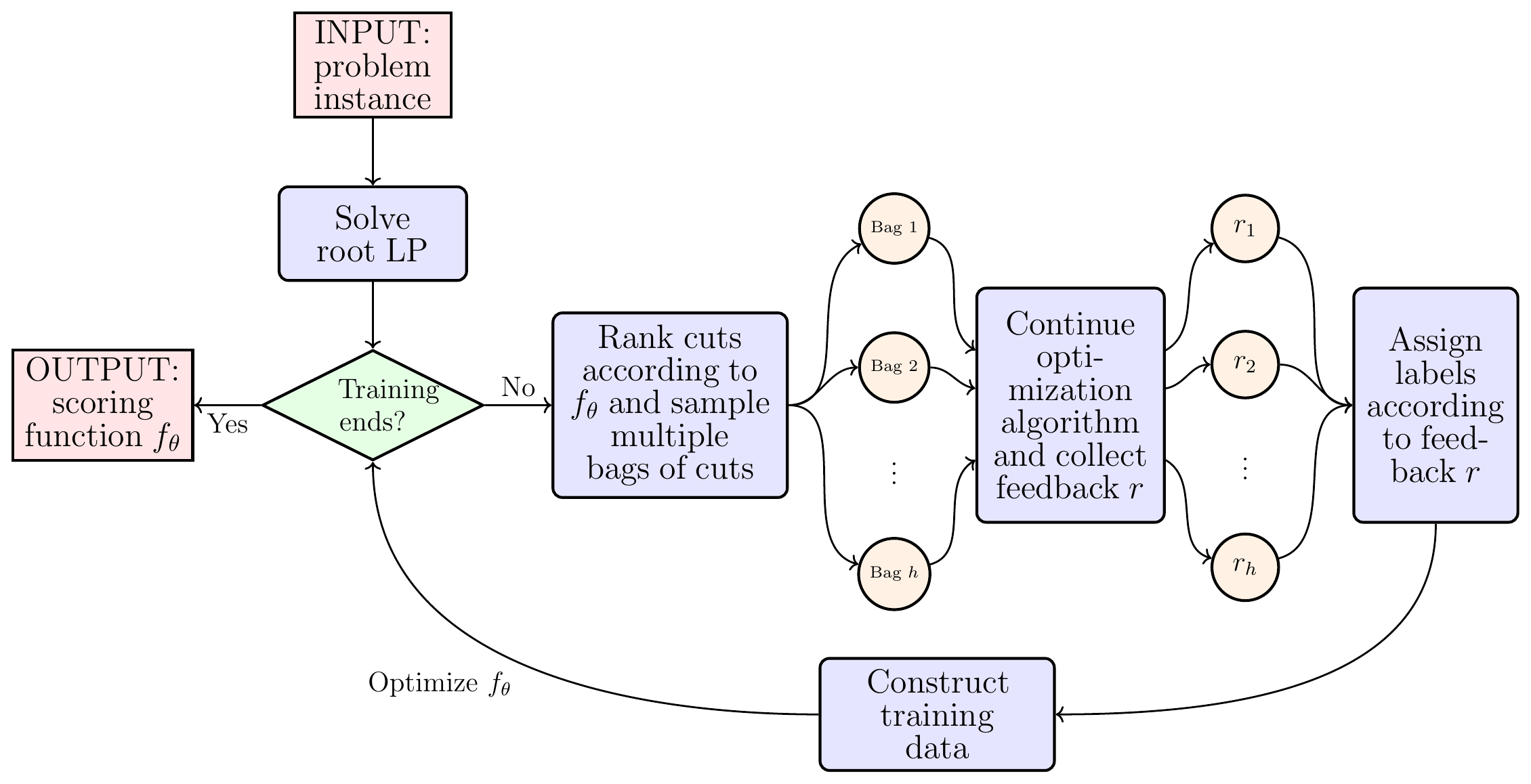}
\caption{The flow chart of active sampling.}
\label{fig:flow_chart}
\end{center}\end{figure}

 \paragraph{\textbf{Active Sampling}}
 In this case, we select the cuts using a pre-trained cut selection policy, which can lead to more promising training samples compared with random sampling. However, collecting samples only based on the pre-trained policy will reduce the sample diversity, which may result in learning a sub-optimal policy. To alleviate this issue, we adopt an $\epsilon$-greedy policy \cite{sutton2018reinforcement}, which is a common approach in reinforcement learning to balance the exploration and exploitation:
 \begin{equation}
 	C' = \left\{
 	\begin{array}{lcl}
 		\text{sample from policy $\pi$}, & &{\text{with probability $1-\epsilon$}} \\
 		\text{sample randomly}, &          &{\text{with probability $\epsilon$}}
 	\end{array}
 	\right.
 \end{equation}
where $C'$ is the selected cut subset, $\pi$ is the cut selection policy derived from the model. During the active sampling phase, the cut selection policy is still being refined using the training samples collected in this phase. The flow chart is displayed in Figure \ref{fig:flow_chart}. Note that active sampling is on-policy, that is, we improve the same policy which is used to collect samples. However, the whole framework of \textsc{Cut Ranking} is in an offline settings, that is, during the test phase, we do not continue to train our model on new MIP instances.

\subsubsection{Constructing Bag Features from Cut Features}
\label{sec:features}
Since the bag features are constructed from the aggregated cut features, thus the first issue is the specification of cut features. To enable better generalization of the model, we design 14 problem-specific atomic features for the cut selection task. Similar to the features of branching variables as provided in Khalil et al. \cite{khalil2016learning}, here we list our designed features for each candidate cut. Specifically, for a given MIP instance with problem property $P=\{x^{*}_{LP}, M\}$, for any generated cut $c_i$: $\alpha_i^\top x \leq \beta_i$, its cut features are shown in Table \ref{tab:cut_features}. 

\begin{table}[!htb]
	\centering
	\caption{Cut features' descriptions and counts.}
	\begin{threeparttable}
	    \small
		\begin{tabular}{l|p{7cm}|c}\toprule
			Feature & Description & Count\\ \midrule
			stats. for cut coeffs. & the mean, max, min, stdev. of cut coefficients $\alpha_i$ & 4\\
			\hline
			stats. for obj. coeffs. & the mean, max, min, stdev. of objective coefficients of cut variables & 4 \\
			\hline
			support &  the proportion of non-zero coefficients in $\alpha_i$ & 1  \\
			\hline
			integral support & the proportion of non-zero coefficients w.r.t integer-constrained variables in $\alpha_i$ & 1 \\
			\hline
			normalized violation & $\max\{0, \frac{\alpha_i^\top x^{*}_{LP}-\beta_i}{|\beta_i|}\}$, which measures the cut violation of the present LP solution & 1\\
			\hline
			distance & $\frac{|\alpha_i^\top x^{*}_{LP} - \beta_i|}{|\alpha_i|}$, which measures the Euclidean distance between the present LP solution and the hyperplane $\alpha_i^\top x = \beta_i$ determined by the cut & 1\\
			\hline 
			parallelism & $\frac{z^\top\alpha_i}{|z|\cdot|\alpha_i|}$, which measures the parallelism between the objective function and the cut & 1\\
			\hline
			expected improvement & $\frac{|z|\cdot|\alpha_i^\top x^{*}_{LP} - \beta_i|}{|\alpha_i|}$, which is an estimation for objective improvement with the cut& 1 \\
			\bottomrule
		\end{tabular}
	\end{threeparttable}
	\label{tab:cut_features}
\end{table}

The top two atomic features correspond to the statistical information of coefficients related to the cut, which help to capture the structural information of the cut. The other features measure the cut characteristics through different measurements as mentioned by Wesselmann and Suhl \cite{cut_selection}, which capture the association between the MIP property and the cut. Moreover, all the above atomic features can be computed quite efficiently, and thus making the time to construct features negligible in the algorithm.

Now that we are able to construct the bag features. For a generated MIP instance, we collect a number of training samples using certain sampling strategies. For each training sample with a selected cut subset $C'$, we first compute and collect the cut features for each cut within $C'$, and obtain the set of corresponding cut features $U_{C'}=\{u_{1}',u_{2}',\ldots, u_{|C'|}'\}$. To prevent feature dimensions from being on different scales, we apply Z-score normalization to re-scale the feature values among all the training samples collected from the same MIP instance. Finally, we introduce a feature mapping function $\phi$, which maps the aggregated cut features to the original cut feature space, and we obtain the final bag features as
\begin{equation}
	x_{C'} = \phi(u_{1}',u_{2}',\ldots, u_{|C'|}').
\end{equation}
The mapping function can be designed to capture the association between cuts within a bag, while in this work, we define $\phi$ as an average function for simplicity, which calculates an average of the cut features. In other works \cite{wang2008aerosol, pappas2014explaining} which are under the collective assumption in MIL settings, the similar average or weighted function is also adopted. Moreover, our empirical studies will show that such a mapping function is effective.

\subsubsection{Assigning Ranking Labels to Training Samples}
\label{sec:label_assign}
The labels for our cut selection task are defined to be binary, taking the value of 1 or 0. The positive labels are assigned to high-quality cut subsets, which are preferable for the selection policy. To measure the quality of a sampled bag, we use the pre-defined reduction ratio of solution time as the feedback signal, and present a labeling scheme based on the ranking of bags. Note that deriving a precise labeling scheme requires collecting all the possible bags for a MIP instance, which is infeasible in general cases. Therefore, we adopt an approximation method, that is to sample a number of bags for each MIP instance, and assign the corresponding bag labels according to their rankings. Empirically, we find that the training is stable under such a labeling scheme.

For training bags $X=\{x_1,x_2,\ldots,x_{h'}\}$ collected from the same MIP instance, let $R=\{r_1, r_2, \ldots, r_{h'}\}$ be the set of corresponding feedbacks. For each sample $x_j$, we define its label as
\begin{equation}
	 y_j = \left\{
	\begin{array}{lcl}
    		1, & \quad\quad {\text{$r_j$ ranks in the $\lambda$\% highest feedbacks}}\\
		0, & \quad\quad{\text{otherwise}}
	\end{array}
	\right.
\end{equation} 
where $\lambda \in (0, 100)$ is a tunable hyper-parameter, which controls the percentage of positive samples. 

The above ranking based binary labeling scheme is suitable for cut selection since the main goal is to distinguish the high-quality cuts between the poor-quality ones without  a full ranking of all candidate cuts.

\subsection{Learning a Scoring Function}
In this subsection, we first introduce the basic architecture of the scoring function, and then present how to train the model.
\subsubsection{Scoring Function Architecture}
\label{function_architecture}
We parameterize the scoring function as a neural network $f_{\theta}(x)$ with parameters $\theta$. In the training phase, for each training bag, the model takes its bag features $x$ as inputs, and outputs a probability distribution with two dimensions, which correspond to $P(y=1 | x)$ and $P(y=0 | x)$, respectively.

The architecture of the model consists of two parts. First, it extracts the bag embeddings $x^e$ through a multi-layer perceptron (MLP), and then feed the embeddings to a softmax layer to output the probability distribution as
\begin{equation}
\begin{aligned}
	x^e &= MLP(x;\theta) \\
	P(y|x) &= \text{softmax}(x^e).
\end{aligned}
\end{equation}

Since the same form of bag features and the cut features, we are able to apply the scoring function at the cut level in the test phase.

\subsubsection{Training the Scoring Function}
We define the probability $P(y=1|x)$ output by the model as the score for the input bag or the cut. The goal is to train the model to output higher scores for positive samples, or high-quality cuts. Given a set of training data $\{(x^j_i,y^j_i)_{i=1}^{h_j}\}_{j=1}^{N}$ collected from $N$ problem instances, we optimize the model parameters $\theta$ to minimize the loss function comprising the cross-entropy loss and the regularization loss as
\begin{equation}
	\label{equ:loss}
	L(\theta) = L_{ce}(\theta) + \gamma \Omega(\theta),
\end{equation}
where $\gamma$ is a hyper-parameter for regularization penalty, and $L_{ce}$ is the cross-entropy loss as
\begin{equation}
	L_{ce}(\theta) = -\sum_{j=1}^{N}\sum_{i=1}^{h_j}y_i^j\log p_{\theta}(y_i^{j}=1|x) + (1-y_i^j)\log p_{\theta}(y_i^{j}=0|x).
\end{equation}
To avoid the overfitting of the model, we adopt a common L2 regularization into the model, which places restrictions on the model parameters as
\begin{equation}
	\Omega(\theta) = \|\theta\|_{2}^{2}.
\end{equation}

Finally, the major steps of our \textsc{Cut Ranking} method are summarized in Algorithm \ref{algo:algo}. The algorithm consists of three phases: data collection phase, training phase and the test phase. In the data collection phase, we collect and construct the training bags and labels; in the training phase, we train our scoring function in a supervised fashion; in the test phase, we apply the scoring function to each candidate cut, and select a cut subset with the highest scores.

\begin{algorithm}[t]
	\caption{\textsc{Cut Ranking} in MIP settings.} 
	\label{algo:algo}
	\small
	\begin{algorithmic}[1]%
		\State \textit{\bf Data Collection Phase:}
		\State \quad Randomly generate a set of training instances $D_{train}$.
		\State \quad For each instance in $D_{train}$:
		\State \quad \quad Sample training bags using strategies mentioned in Section \ref{sec:strategy};
		\State \quad \quad Construct the bag features via Section \ref{sec:features};
		\State \quad \quad Assign labels to the bags via Section \ref{sec:label_assign}.
		\State \textit{\bf Training Phase:}
		\State \quad Initialize the scoring function $f_\theta$ and the learning rate $\mu$.
		\State \quad Repeat 
		\State \quad \quad For each batch in training data:
		\State \quad \quad \quad Calculate the loss $L(\theta)$ in Equation \ref{equ:loss} ;
		\State \quad \quad \quad Optimize model parameters: $\theta \leftarrow \theta - \mu \nabla_{\theta}L(\theta)$;
		\State \quad until $\theta$ converge.
		\State \textit{\bf Test Phase:}
		\State \quad For each instance in the test set $D_{test}$:
		\State \quad \quad Obtain the scores for each candidate cut using $f_{\theta}$;
		\State \quad \quad Select the top $K$\% cuts with the highest scores.
	\end{algorithmic}
\end{algorithm}

\section{Experiments}
\subsection{Experimental Setup}
\subsubsection{Benchmarks}
Our benchmarks consist of synthetic MIP problems and the real-world production planning problems. The synthetic MIPs consist of four classical classes of problems: Set Cover, Knapsack, Planning and General MIP. For the ease of data collection, we find the problem instances which are solvable within 25 seconds. For each class of problems, we randomly generate 100 training instances and 30 test instances using different random seeds. 

The large-scale real-world daily production planning problems are divided into two phases, the \textsl{Offline Phase} during January 2021, and the \textsl{Online Phase} during March 2021. These two phases correspond to the offline datasets and online datasets for our experiments.

\subsubsection{Metrics}
The cut selection policy will affect the branch-and-cut algorithm from two aspects, the solution time and the number of nodes visited. Compared with the cases without cuts, we define two metrics for cut selection based on the above algorithm feedbacks: \textbf{the reduction ratio of solution time}; and \textbf{the reduction ratio of the number of nodes visited}. 

Notably, the solution time is not directly determined by the number of nodes visited, since we need to also consider the solving time of each node relaxation. Therefore, the primary metric to qualify the cut selection policy is the reduction ratio of solution time. For each conducted experiment, we show the mean and standard deviation of results on the test set, and highlight the best average results. 

\subsubsection{Baselines}
For synthetic datasets, we examine our \textsc{Cut Ranking} module against five widely used manually-designed heuristics for cut selection, including:
\begin{itemize}
	\item \textsc{Random}: select cuts according to a stochastic policy.
	\item \textsc{Violation}: select the cuts with larger violation.
	\item \textsc{Normalized Violation}: select the cuts with larger normalized violation.
	\item \textsc{Distance}: select the cuts with larger Euclidean distance from the root LP solution $x^{*}_{LP}$.
	\item \textsc{Parallelism}: select the cuts which are more parallel to the objective function.
\end{itemize}

For real-world datasets, we compare \textsc{Cut Ranking} with the fine-tuned manually heuristics which are adopted in Huawei's proprietary solver.

Note that RL2C~\cite{tang2020reinforcement} is not included in the baselines since their formulation is based on sequential decision making, and the main goal of RL2C is to reduce the total number of added cuts to solve the MIP to optimality. The main algorithmic framework of RL2C is the cutting plane method, which iteratively adds cuts to the initial LP relaxation. Specifically, in RL2C, a new cut is selected and added from the candidate cut set at each step and the solver is required to execute the LP instantly to obtain the change of objective value $\mathbf{z}^\top \mathbf{x}$. However, our defined task is essentially a one-step decision problem, and the main goal of \textsc{Cut Ranking} is to improve the solution time of MIP in the algorithmic framework of branch-and-cut. Therefore, RL2C is highly incompatible with our settings, and its RL formulation is impractical for large-scale MIPs since it may lead to much more sampling and computational costs. Moreover, the policy architecture of RL2C is based on LSTM, and the network inputs include all the constraints and available candidate cuts, which also makes it infeasible to apply in large-scale scenarios.

\subsubsection{Implementation of Algorithms}
The algorithmic framework for synthetic MIPs is branch-and-cut. We implement the vanilla branch-and-bound algorithm, and use the open-source solver Python-MIP \cite{saltzman2002coin} for solving LPs and generating cuts. For the real-world datasets, the optimization is based on a proprietary industrial solver of Huawei Company. We enable cut selection before the rounding procedure.

\subsubsection{Hyper-parameters}
\paragraph{Policy Architecture}
The implemented policy network of \textsc{Cut Ranking} is a 4-layer fully-connected neural network, including an input layer, two hidden layers with 30 and 15 hidden units respectively and tanh activation, and an output layer. As mentioned in Section \ref{function_architecture}, the input layer accepts the cut (or bag) features with dimension 14 as the network inputs. The output layer outputs a probability distribution with two dimensions, and we define the positive probability as the score for the input cut (or bag).
\paragraph{Hyper-parameters of Algorithms}

For both synthetic datasets and real-world datasets, we set the hyper-parameter $K$ to 30, $\lambda$ to 50, and $\gamma$ to $0.1$ after hyper-parameter tuning. For synthetic datasets, the number of sampled bags for each MIP instance is 100, and the total number of collected training samples is $10^4$. For real-world datasets, to improve the sample efficiency, the exploring policy starts from a well-tuned heuristic combined with $\epsilon$-greedy. The training samples are collected from the daily production planning problems within a month.

\subsection{Experiments on Synthetic MIP Datasets}
\subsubsection{Experiment \uppercase\expandafter{\romannumeral1}: the Quality of Selected Cuts}
To check the effectiveness of our proposed ranking-based cut selection policy, we conduct comparative experiments on four classes of MIP problems. For Set Cover, the number of elements and sets are both set to 200, and the resulting problem size is $200 \times 200$; for Knapsack, the number of items is 700, and the resulting problem size is $701 \times 700$; for Planning, the number of factories and demands is 20 and 50, respectively, and the resulting problem size is $140 \times 1420$; for General MIP, owing to its complex problem structure, the problem size is set to be $30\times 30$. Although the problem size varies from class to class, the mean solution time (without cuts) is close, thus the difficulty of MIP instances for each class is at the same level. For each class of instances, the number of generated candidate cuts is roughly 20.

As shown in Table \ref{tab:exp_1}, in terms of the problem solving time, our proposed \textsc{Cut Ranking} policy has achieved higher average reduction ratio of solution time over other baseline policies on all the problems, which leads to less solution time. Moreover, the \textsc{Cut Ranking} policy has also shown to be more stable on multiple instances since the standard deviation is relatively smaller compared to the mean. For the Knapsack and Planning problems, the average performance of the heuristic \textsc{Normalized Violation} is comparable to us, while with a large variance, which indicates that such a heuristic often suffers from performance fluctuations. The results are similar for other human-designed heuristics, compared to the \textsc{Cut Ranking} policy, they have shown larger performance variance, and may slow down the solving process in many cases. 

Considering the impact of cut selection on the size of the branch-and-cut search tree, for the \textsc{Cut Ranking} policy, the number of nodes visited has decreased more significantly than other baselines over the Set Cover and General MIP problem instances. For the Knapsack and Planning problems, our proposed policy has also shown to be competitive, and also achieve smaller variance.

Overall, these results indicate that the \textsc{Cut Ranking} policy has improved the optimization algorithm more significantly compared to the human-designed heuristics. Moreover, the \textsc{Cut Ranking} policy is capable of speeding up the solving time for all the generated test cases, while the effects of other cut selection heuristics suffer from instability.

\begin{table}[htb]  
	\centering  
    \caption{Evaluation results of cut selection policies in terms of the reduction ratio of solving time (higher is better), and the reduction ratio of visited nodes (higher is better).}
    \label{tab:exp_1}  
	\resizebox{1.0\textwidth}{!}{
	\begin{threeparttable} 
		\begin{tabular}{ccccccccc}
			\toprule  
			\multirow{2}{*}{Method}&  
			\multicolumn{2}{c}{Set Cover}&\multicolumn{2}{c}{ Knapsack}&\multicolumn{2}{c}{Planning}&\multicolumn{2}{c}{General MIP}\cr  
			\cmidrule(lr){2-3} \cmidrule(lr){4-5} \cmidrule(lr){6-7} \cmidrule(lr){8-9}
			 & Time & Nodes & Time & Nodes & Time & Nodes & Time & Nodes\cr  
			\midrule  
			\textsc{Random} & 0.09$\pm$0.45 & 0.27$\pm$0.48 &  0.14$\pm$0.27 & 0.62$\pm$0.67 & 0.07$\pm$0.14 & 0.40$\pm$0.30&-0.10$\pm$0.11 & -0.03$\pm$0.12\cr  
			\textsc{Violation} & 0.17$\pm$0.19 & 0.35$\pm$0.18 & 0.21$\pm$0.33  & 0.50$\pm$0.88 &  -0.01$\pm$0.24& 0.36$\pm$0.35& 0.11$\pm$0.19& 0.25$\pm$0.20\cr  
			\textsc{Norm-Violation} &0.16$\pm$0.22  & 0.34$\pm$0.22  & 0.25$\pm$0.24 & \textbf{0.70}$\pm$0.38 & 0.17$\pm$0.20 & \textbf{0.48}$\pm$0.30 & 0.13$\pm$0.39 & 0.20$\pm$0.35\cr  
			\textsc{Distance} & 0.12$\pm$0.36 & 0.30$\pm$0.40 & 0.10$\pm$0.19 & 0.41$\pm$0.44&  0.01$\pm$0.15&0.35$\pm$0.33 & 0.10$\pm$0.28 &0.25$\pm$0.20 \cr  
			\textsc{Parallelism} & 0.06$\pm$0.34 & 0.23$\pm$0.37& 0.16$\pm$0.18 & 0.58$\pm$0.34 & -0.06$\pm$0.17 &0.34$\pm$0.36& 0.03$\pm$0.20& 0.26$\pm$0.20\cr  
			\textsc{Cut Ranking} & \textbf{0.21}$\pm$0.16 & \textbf{0.49}$\pm$0.16 & \textbf{0.27}$\pm$0.20 & 0.69$\pm$0.32 & \textbf{0.18}$\pm$0.17 & \textbf{0.48}$\pm$0.28 & \textbf{0.32}$\pm$0.25& \textbf{0.38}$\pm$0.20 \cr 
			\bottomrule  
		\end{tabular}
	\end{threeparttable}}  
\end{table}

\subsubsection{Experiment \uppercase\expandafter{\romannumeral2}: Study of Generalization Ability}
To test if our proposed policy has the generalization ability over problems with different structures or scales, we conduct three experiments to try to answer the following questions:
\begin{itemize}
	\item Can the \textsc{Cut Ranking} policy generalize to the same class of problems with different sizes?
	\item Can the \textsc{Cut Ranking} policy generalize to the same class of problems with different coefficient ranges?
	\item Can the \textsc{Cut Ranking} policy generalize to the problems with different structures?
	
\end{itemize}
\paragraph{Problem Size}
Table \ref{tab:exp_gen_scale} presents the results on Knapsack problems with different problem sizes. The learning-based policy is trained on the problem instances with 700 items, and tested against other heuristics on instances with 600, 800, 900 and 1000 items. As can be seen from the table, the \textsc{Cut Ranking} policy has shown a higher reduction ratio of solution time over other baselines on test instances with 600, 800 and 1000 items. For problem instances with 800 items, although the heuristic \textsc{Violation} and \textsc{Parallelism} have achieved a slightly higher averaged reduction ratio, our policy has a much lower variance, thus is more preferable for cut selection. 

The \textsc{Cut Ranking} policy also results in a smaller branching tree compared to most baseline heuristics. From these two aspects, we conclude that our proposed policy trained on problem instances of a certain scale can be applied to the same class of problems of different scales.
\begin{table}[htb]  
	\centering  
	\caption{Evaluation results of cut selection policies on Knapsack problems with different scales. We train our learning module on 100 randomly generated Knapsack instances with 700 items.}
	\label{tab:exp_gen_scale} 
	\resizebox{1.0\textwidth}{!}{
		\begin{threeparttable}  
			\begin{tabular}{ccccccccc}  
				\toprule  
				\multirow{2}{*}{Method}&  
				\multicolumn{2}{c}{600 items}&\multicolumn{2}{c}{800 items}&\multicolumn{2}{c}{900 items}&\multicolumn{2}{c}{1000 items}\cr  
				\cmidrule(lr){2-3} \cmidrule(lr){4-5} \cmidrule(lr){6-7} \cmidrule(lr){8-9}
				& Time & Nodes & Time & Nodes & Time & Nodes & Time & Nodes\cr  
				\midrule  
				\textsc{Random} & 0.20$\pm$0.29  & 0.69$\pm$0.28& 0.20$\pm$0.24  & 0.64$\pm$0.31 &0.26$\pm$0.41& 0.65$\pm$0.70& 0.21$\pm$0.20 & 0.55$\pm$0.40\cr  
				\textsc{Violation}& 0.27$\pm$0.31  & 0.75$\pm$0.30 & 0.21$\pm$0.25 & 0.67$\pm$0.32 & 0.31$\pm$0.27& \textbf{0.77}$\pm$0.27 & 0.23$\pm$0.24& 0.50$\pm$0.54\cr  
				\textsc{Norm-Violation} & 0.18$\pm$0.31  & \textbf{0.75}$\pm$0.31  & 0.30$\pm$0.26 & 0.74$\pm$0.33 & 0.24$\pm$0.28& 0.71$\pm$0.27& 0.17$\pm$0.24  & 0.50$\pm$0.53\cr  
				\textsc{Distance} & 0.21$\pm$0.28 & 0.62$\pm$0.35 & 0.29$\pm$0.27 & 0.67$\pm$0.39 &0.22$\pm$0.39 & 0.65$\pm$0.60& 0.22$\pm$0.22 & 0.57$\pm$0.45 \cr  
				\textsc{Parallelism} & 0.10$\pm$0.24 & 0.58$\pm$0.40&  0.27$\pm$0.26& 0.66$\pm$0.38 & \textbf{0.33}$\pm$0.40 & 0.74$\pm$0.68& 0.19$\pm$0.20 & 0.53$\pm$0.41 \cr  
				\textsc{Cut Ranking} & \textbf{0.28}$\pm$0.15 & 0.74$\pm$0.22& \textbf{0.33}$\pm$0.24 & \textbf{0.80}$\pm$0.31 & 0.31$\pm$0.16 &0.71$\pm$0.18 & \textbf{0.25}$\pm$0.18 & \textbf{0.61}$\pm$0.38 \cr 
				\bottomrule  
			\end{tabular}  
	\end{threeparttable}} 
\end{table} 

\paragraph{Coefficient Ranges}
The parameters of Knapsack problems include \textsl{the maximal number}, \textsl{the maximal value} and \textsl{the maximal weight} for one type of item, which restricts the range of each randomly generated coefficient. We generate four sets of Knapsack instances with parameters set to 10, 20, 50 and 100. We train our cut selection module on the instances with parameters set to 10, and test on other three sets of instances.

The results are set out in Table \ref{tab:exp_gen_coeff}, from which we can observe that the \textsc{Cut Ranking} policy clearly outperforms other baselines on problem instances with coefficients range between 0 and 20. For problems with coefficients range between 0 and 50, the \textsc{Cut Ranking} policy is still superior to the baselines, while the performance gap between them has decreased much. Such a phenomenon is more striking for problems with larger ranges of coefficients, as demonstrated in the rightmost two columns of Table \ref{tab:exp_gen_coeff}, the \textsc{Cut Ranking} policy does not have clear advantages over the heuristic \textsc{Normalized Violation}.
 
Taken together, the results reveal that our \textsc{Cut Ranking} policy can generalize to the problems with different coefficient ranges. However, problems with a large range of coefficients will limit the generalization ability of the learned policy.

\begin{table}[htb]  
	\centering  
	\caption{Evaluation results of cut selection policies on Knapsack problems with the same size (700 items) but different ranges of coefficients. We train our learning module on Knapsack instances with coefficients range between $(0,10]$.}
	\label{tab:exp_gen_coeff}  
	\resizebox{1.0\textwidth}{!}{
		\begin{threeparttable}  
			\begin{tabular}{ccccccc}  
				\toprule  
				\multirow{2}{*}{Method}&  
				\multicolumn{2}{c}{$0<\text{coeff.}\leq 20 $}&\multicolumn{2}{c}{$0<\text{coeff.}\leq 50$}&\multicolumn{2}{c}{$0< \text{coeff.} \leq 100$}\cr  
				\cmidrule(lr){2-3} \cmidrule(lr){4-5} \cmidrule(lr){6-7} 
				& Time & Nodes & Time & Nodes & Time & Nodes\cr  
				\midrule  
				\textsc{Random} & 0.42$\pm$0.54 & 0.63$\pm$0.87 & 0.36$\pm$0.35 & 0.70$\pm$0.58 & 0.14$\pm$0.27 & 0.62$\pm$0.67\cr  
				\textsc{Violation}& 0.43$\pm$0.51 & 0.64$\pm$0.80  & 0.35$\pm$0.40 & \textbf{0.71}$\pm$0.44  & 0.21$\pm$0.33 & 0.50$\pm$0.88 \cr  
				\textsc{Norm-Violation} & 0.41$\pm$0.52 & 0.62$\pm$0.80  & 0.37$\pm$0.27 & 0.67$\pm$0.40 & \textbf{0.25}$\pm$0.24 & \textbf{0.70}$\pm$0.38\cr  
				\textsc{Distance} & 0.30$\pm$0.32  & 0.56$\pm$0.38 &  0.32$\pm$0.26& 0.68$\pm$0.30 & 0.10$\pm$0.19 & 0.41$\pm$0.44  \cr  
				\textsc{Parallelism} & 0.33$\pm$0.34& 0.58$\pm$0.41 & 0.33$\pm$0.28& 0.65$\pm$0.34 & 0.16$\pm$0.18 & 0.58$\pm$0.34   \cr  
				\textsc{Cut Ranking} & \textbf{0.52}$\pm$0.30 & \textbf{0.86}$\pm$0.23& \textbf{0.39}$\pm$0.27& 0.68$\pm$0.39 & 0.23$\pm$0.19 & 0.65$\pm$0.30 \cr 
				\bottomrule  
			\end{tabular}  
	\end{threeparttable}} 
\end{table}

\paragraph{Problem Structure}
To explore the generalization ability of the \textsc{Cut Ranking} policy on problems with different structures, we test the policy trained on Knapsack instances on Set Cover, Planning and General MIP problems. The generated problem size in this experiment is kept the same as experiment \uppercase\expandafter{\romannumeral1}. 

From the results shown in Table \ref{tab:exp_gen_class}, the \textsc{Cut Ranking} policy outperforms other baselines on Set Cover problems, which shows greater average improvement with a lower variance. For Planning problems, the heuristic \textsc{Normalized Violation} has achieved slightly better average performance compared to the learned policy, however, the variance is much larger. For General MIPs, the \textsc{Cut Ranking} policy also outperforms other baselines, nevertheless, it fails to reduce the solving time of some test instances. In view of the more complex problem structure for General MIPs, it is more difficult for the learned policy to generalize to these problem instances.

In summary, these results show that the \textsc{Cut Ranking} policy has certain generalization ability on problems with different structures. The results also indicate that the General MIPs and the Knapsack problems may have a relatively greater distinction in problem structures, which may lead to increased divergence between the distribution of training and test dataset.

\begin{table}[htb]  
	\centering  
	\caption{Evaluation results of cut selection policies on different classes of problems. We trained our learning module on Knapsack instances, and test on other classes of problem instances which at the same difficulty level.}
	\label{tab:exp_gen_class}  
	\resizebox{1.0\textwidth}{!}{
	\begin{threeparttable}
		\begin{tabular}{ccccccc}  
			\toprule  
			\multirow{2}{*}{Method}&  
			\multicolumn{2}{c}{Set Cover}&\multicolumn{2}{c}{Planning}&\multicolumn{2}{c}{General MIP}\cr  
			\cmidrule(lr){2-3} \cmidrule(lr){4-5} \cmidrule(lr){6-7} 
			& Time & Nodes & Time & Nodes & Time & Nodes\cr  
			\midrule  
			\textsc{Random} & 0.09$\pm$0.45 &0.27$\pm$0.48  &0.07$\pm$0.14  & 0.40$\pm$0.30 & -0.10$\pm$0.11  & -0.03$\pm$0.12 \cr  
			\textsc{Violation}& 0.17$\pm$0.19& 0.35$\pm$0.18 & -0.01$\pm$0.24 & 0.36$\pm$0.35 & 0.11$\pm$0.19 &0.25$\pm$0.20 \cr  
			\textsc{Norm-Violation} & 0.16$\pm$0.22 & 0.34$\pm$0.22  & \textbf{0.17}$\pm$0.20  & \textbf{0.48}$\pm$0.30 & 0.13$\pm$0.39 & 0.20$\pm$0.35 \cr  
			\textsc{Distance} & 0.12$\pm$0.36 & 0.30$\pm$0.40 & 0.01$\pm$0.15 & 0.35$\pm$0.33 & 0.10$\pm$0.28 & 0.25$\pm$0.20 \cr  
			\textsc{Parallelism} & 0.06$\pm$0.34 & 0.23$\pm$0.37& -0.06$\pm$0.17 & 0.34$\pm$0.36 & 0.03$\pm$0.20 & \textbf{0.26}$\pm$0.20  \cr  
			\textsc{Cut Ranking} & \textbf{0.22}$\pm$0.18  & \textbf{0.52}$\pm$0.17 & 0.16$\pm$0.11 & 0.33$\pm$0.26 & \textbf{0.15}$\pm$0.23 & 0.23$\pm$0.18 \cr 
			\bottomrule  
		\end{tabular}  
	\end{threeparttable}}
\end{table}

\subsection{Experiments on Large-Scale Real-World Tasks}
To further evaluate the quality of the proposed \textsc{Cut Ranking} policy, we embed the cut selection module in an industrial large-scale MIP solver developed by Huawei Company, and conduct both offline and online experiments on the real-world production planning problems with more than $10^7$ variables and constraints daily. To our knowledge, this is the first study to apply machine learning into cut selection for large-scale MIPs with more than $10^7$ variables and constraints.

\paragraph{MIP Statistics}
The daily production-planning problems have similar problem structures, while problem properties change day by day. For a better description of problem characteristics, we record the number of variables and constraints, the density of the constraint matrix, the mean and standard deviation of objective coefficients, and the number of integer variables. 

\begin{figure}[t]\begin{center}
\includegraphics[width=1.0\textwidth]{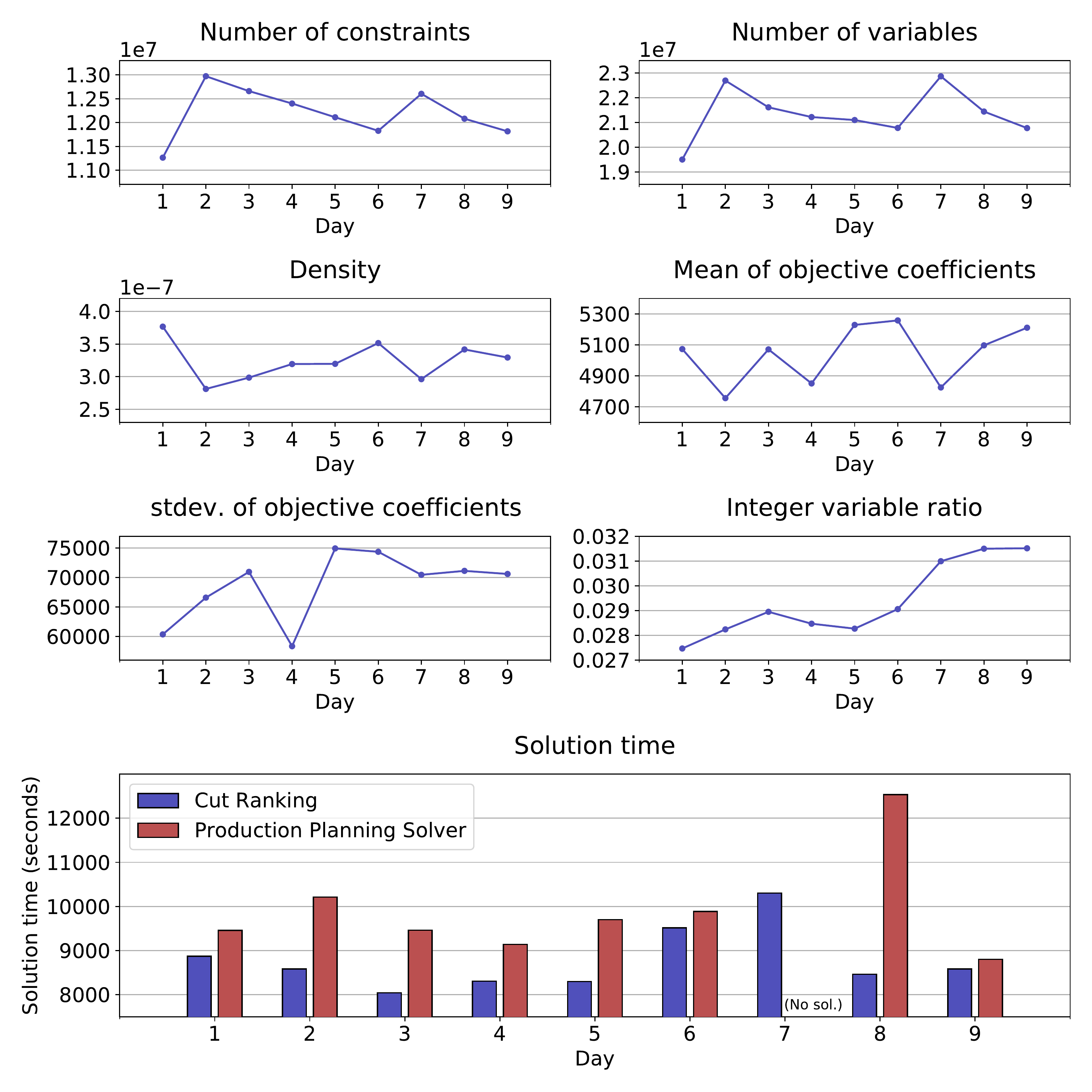}
\caption{The daily MIP statistics and the evaluation results for offline datasets.}
\label{fig:offline_stats}
\end{center}\end{figure}

\begin{figure}[t]\begin{center}
\includegraphics[width=1.0\textwidth]{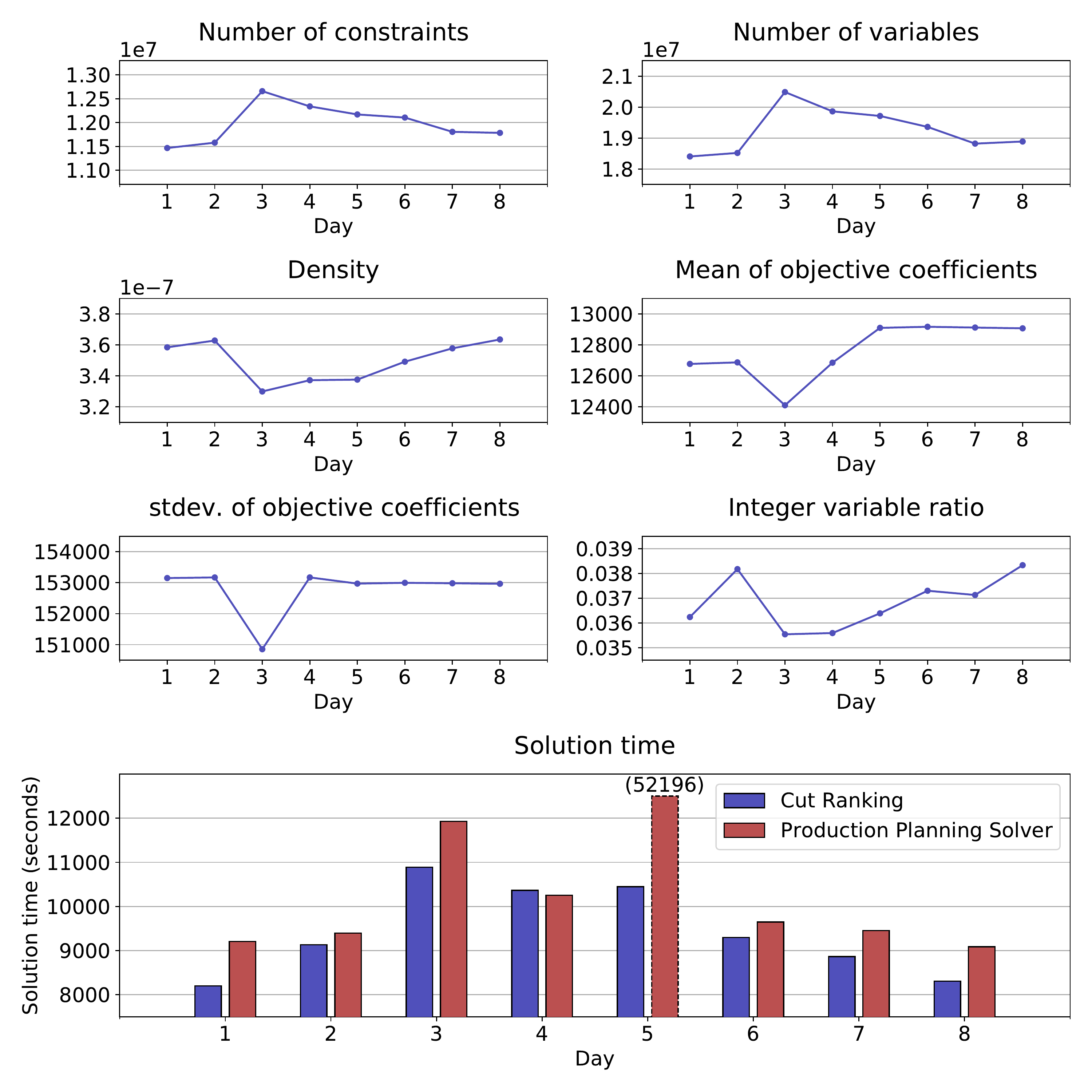}
\caption{The daily MIP statistics and the evaluation results for online datasets.}
\label{fig:online_stats}
\end{center}\end{figure}

Table \ref{tab:mean_property} shows the mean MIP statistics of production planning datasets. The upper six sub-figures of Figure \ref{fig:offline_stats} and \ref{fig:online_stats} demonstrate the visualization of problem statistics during January 2021 and March 2021 for offline and online datasets, respectively. Moreover, the number of generated candidate cuts for each problem instance is between 1200 and 1400.

\begin{table}[!htbp]
    \centering
    \caption{The mean MIP statistics of the large-scale production planning tasks.}
     \label{tab:mean_property}
    \resizebox{1.0\textwidth}{!}{
    \begin{threeparttable}
    \begin{tabular}{c|c|c|c|c|c|c} \toprule
    Dataset & \# constraints & \# variables  & density & mean of obj. coeff. & stdev. of obj. coeff. & integer variable ratio \\ \midrule
    Offline & 12,192,747 & 21,334,700 & 3.2$\times 10^{-7}$ & 5,041 &68,642 & 0.030\\    
    \hline
    Online & 11,988,209 & 19,256,236 & 3.5$\times 10^{-7}$& 12,763 &152,782 & 0.037 \\
    \bottomrule 
    \end{tabular}   
    \end{threeparttable}}
   
\end{table}

\paragraph{Comparison Experiments} Our learning module is trained on the offline collected samples, and applied to both offline and online datasets. We compare our \textsc{Cut Ranking} policy against the production planning solver with a manually heuristic for cut selection. As shown in Figure \ref{fig:offline_stats} and \ref{fig:online_stats}, our \textsc{Cut Ranking} policy has led to less solution time on both offline and online datasets without the accuracy loss of solution, and the average speedup ratio has reached 14.98\% and 12.42\%, respectively. Interestingly, for the problem instance of day $7$ in offline datasets, the solver with the manually heuristic is unable to return a solution within a limited time; for the problem instance of day $5$ in online datasets, the solver with the manually heuristic costs much more time to obtain the solution compared to the solver with the learning module. These results have further demonstrated the importance of deriving a proper cut selection policy. Our proposed ranking-based cut selection policy has shown to be more robust and efficient compared to the baseline, and is also generalizable to the daily problems with different properties.

\begin{figure}[t] \centering    
	\subfigure[Parallelism] { 
		\includegraphics[width=0.46\columnwidth]{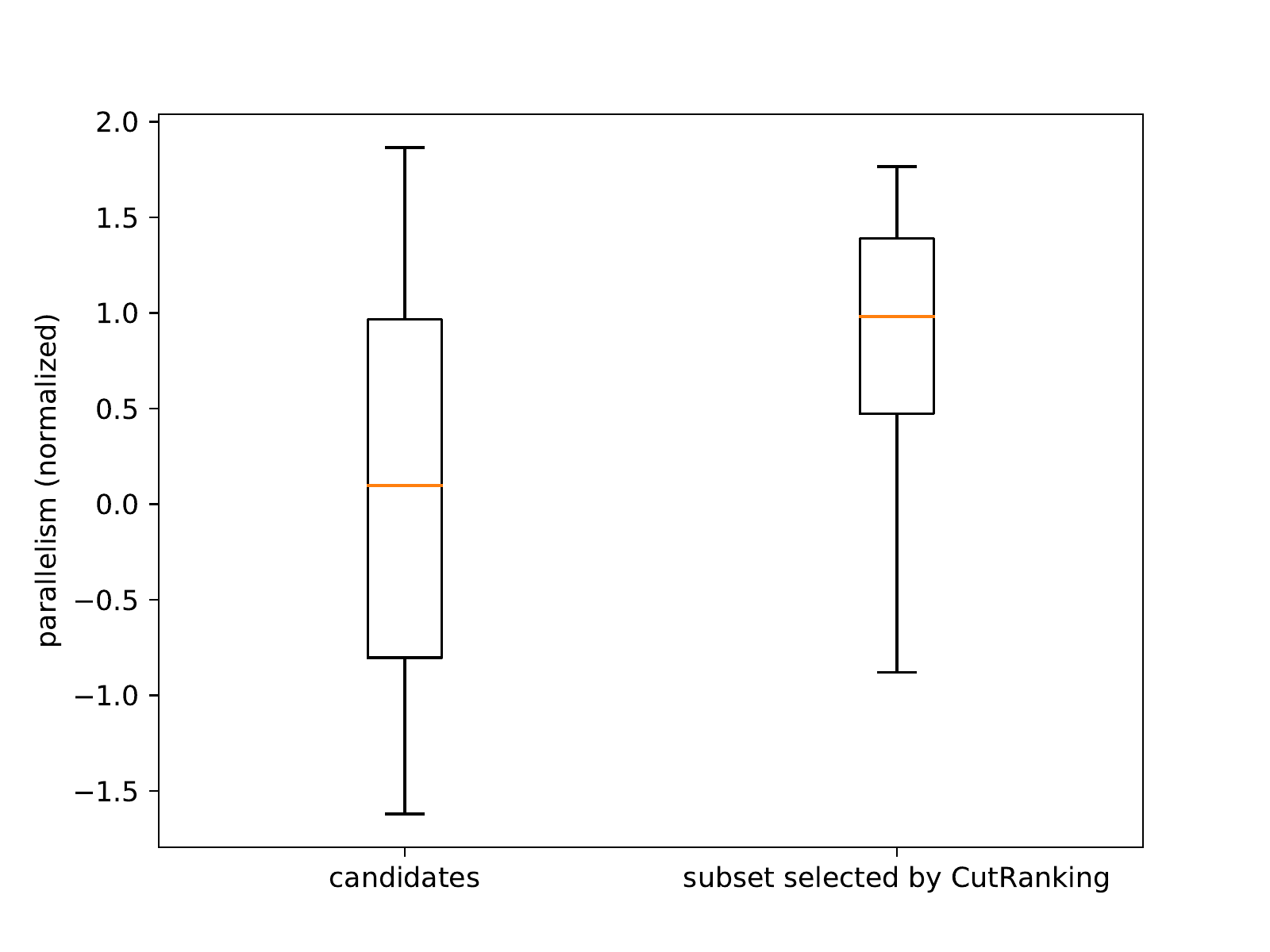}  
	}     
	\subfigure[Stats. for obj. coeff. (mean)] { 
		\includegraphics[width=0.46\columnwidth]{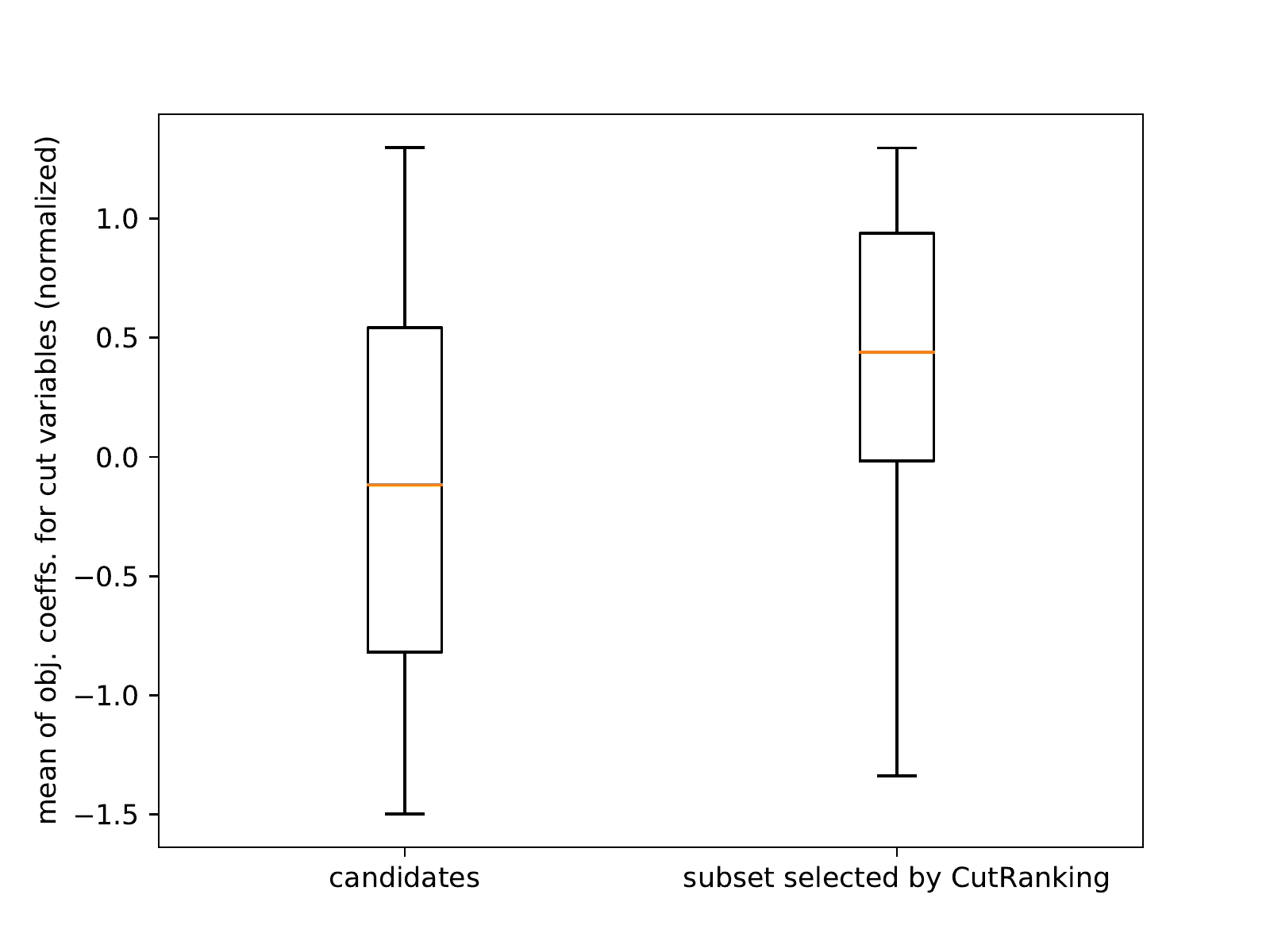}     
	}    
	\caption{Box plot of informative cut features for production planning problems}     
	\label{fig:boxplot}     
\end{figure}

\paragraph{Strategy Captured by \textsc{Cut Ranking}}
\textsc{Cut Ranking} will learn to find the informative cut features or the informative feature combinations for the given class of MIP problems. We analyze \textsc{Cut Ranking} on the real-world production planning problems since they have particular problem structures. As we have tested, the heuristic based on \textsc{Parallelism} performs well on the problem instances, which indicates that \textsc{Parallelism} is one of  the informative cut features for the production planning problems. We find that \textsc{Cut Ranking} has learned to select cuts with larger \textsc{Parallelism} as well. Besides, \textsc{Cut Ranking} also prefers cuts with larger mean value of objective coefficients of cut variables. For a better demonstration, we visualize the informative features found by \textsc{Cut Ranking} using the box plot, which can be seen in Figure \ref{fig:boxplot}.

\section{Conclusion}
In this paper, we presented $\textsc{Cut Ranking}$ method for cut selection in the context of branch-and-cut algorithm for MIPs. To tackle the infeasibility of acquiring the supervised label for a single cut the task, we proposed to model the learning process in the settings of multiple instance learning. Moreover, we designed several problem-specific features for cuts, and provided a scheme to construct bag features and labels for training. The experimental results on synthetic MIPs have demonstrated that our learned ranking-based cut selection policy is more competitive compared to other manually heuristics, and also with generalization ability on problems with different scales, coefficient ranges or structures. For the real-world production planning tasks, our $\textsc{Cut Ranking}$ method has also significantly improved the efficiency Huawei's industrial MIP solver without the accuracy loss of solution, achieving a speedup ratio of 14.98\% and 12.42\% in offline and online A/B testings, respectively. The empirical findings of this work provide a deeper insight into the generalization ability of machine learning techniques for cut selection, and reveal that the machine learning module can be incorporated in the solver to improve the solution process even for large-scale MIPs.

\section*{Acknowledgements} 
Weinan Zhang is supported by “New Generation of AI 2030” Major Project (2018AAA0100900) and National Natural Science Foundation of China (62076161). The work is also sponsored by Huawei Innovation Research Program.

\bibliography{mybibfile}

\end{document}